\documentclass[11pt,draft,a4paper]{article}
\usepackage{amsthm,amstext,amsmath,amscd,amssymb,eucal,latexsym}
\usepackage[matrix,arrow]{xy}

\title
{Irregular elliptic surfaces of degree $12$ \\
in the projective fourspace}
\author{Hirotachi Abo \\
Fachrichtung Matematik \\
Universit\"at des Saarlandes \\
D-66041 Saarbr\"ucken \\
Germany \\
{\it E-mail}: {\tt abo@math.uni-sb.de}
\and
Kristian Ranestad \\
Matematisk Institutt \\
Universitetet i Oslo \\
P.b. 1053 Blindern \\
N-0316 Oslo 3, Norway \\
{\it E-mail}: {\tt ranestad@math.uio.no}
}

\date{}
    \newtheorem{thm}{Theorem}[section]
    \newtheorem{lem}[thm]{Lemma}
    
    \newtheorem{prop}[thm]{Proposition}
\theoremstyle{definition}
    \newtheorem{df}[thm]{Definition}
    \newtheorem{ex}[thm]{Example}
    \newtheorem{rem}[thm]{Remark}

\def\C      {\mathbb{C}}                      

\def\Z      {\mathbb{Z}}

\def\H      {\mathbb{H}}                

\def\P      {\mathbb{P}}                

\def\AA  {\mathcal A}
\def\EE  {\mathcal E}
\def\KK  {\mathcal K}
\def\CC  {\mathcal C}

\def\FF  {\mathcal F}
\def\GG  {\mathcal G}

\def\BB  {\mathcal B}
\def\II  {\mathcal I}
\def\OO  {\mathcal O}

\def\QQ  {\mathcal Q}

\def\H   {\mathrm  H}
\def\h   {\mathrm h}

\def\ker{\mathop{\mathrm{Ker}}\nolimits}
\def\im{\mathop{\mathrm{Im}}\nolimits}
\def\coker{\mathop{\mathrm{Coker}}\nolimits}
\def\rank{\mathop{\mathrm{rank}}\nolimits}

\def\Ext{\mathop{\mathrm{Ext}}\nolimits}

\def\Syz{\mathop{\mathrm{Syz}}\nolimits}
\def\Hom{\mathop{\mathrm{Hom}}\nolimits}

\numberwithin{equation}{section}

\begin{document}

\maketitle

\section{Introduction}
The purpose of this paper is to give
two different constructions of
irregular surfaces of degree~$12$ in~$\P^4$.
These surfaces are new and
of interest in the classification of
smooth non-general type surfaces in ~$\P^4$.
This classification problem is motivated by the theorem
of Ellingsrud and Peskine~\cite{ep} which says that
the degree of smooth non-general type surfaces
in~$\P^4$ is bounded.
Moreover the new family is one of a few known
families of irregular surfaces in~$\P^4$.
The other known irregular smooth surfaces in~$\P^4$
are the elliptic quintic scrolls,
the elliptic conic bundles~\cite{ads},
the bielliptic surfaces of degree~$10$~\cite{serrano}
and degree~$15$~\cite{adhpr},
the minimal abelian surfaces of degree~$10$~\cite{hm}
and the non-minimal abelian surfaces of degree~$15$~(cf.~\cite{aure})
up to pull-backs of these families by suitable
finite maps from~$\P^4$ to~$\P^4$ itself.

The first author came across the irregular elliptic surfaces
when studying a stable rank three vector bundle~$\EE$
on~$\P^4$ with Chern classes~$(c_1,c_2,c_3) = (5,12,12)$
(see~\cite{abo} for the explicit construction of this bundle).
The dependency locus of two sections of~$\EE$
is a smooth surface of the desired type.

Our first construction uses monads.
The basic idea of monads is
to represent a given coherent sheaf
as a cohomology sheaf of a complex of simpler vector bundles.
A useful way to construct monads is
Horrocks' technique of killing cohomology.
As we will see in $\S$~\ref{monad},
for the ideal sheaf of a smooth surface~$X$ in~$\P^4$
this technique is closely related to
the graded finite length modules
$\H^i_* \II_X = \bigoplus_{m \in \Z} \H^i (\P^4, \II_X (m))$, $i = 1,2$,
over the polynomial ring~$\C[x_0, \dots ,x_4]$ called
the {\it Hartshorne-Rao modules} of~$X$.
To construct an irregular elliptic surface,
we construct the ideal sheaf via its Hartshorne-Rao modules.

The second construction uses liaison,
and reduces the construction of an irregular surface as above
to the construction of a simpler
locally complete intersection surface.
More precisely we show that
the elliptic irregular surface of the first construction
is linked~$(5,5)$ to a reducible surface of degree~$13$.
This reducible surface consists of
a singular rational cubic scroll surface in a hyperplane
of~$\P^4$
and a smooth general type surface
of degree~$10$ and sectional genus~$10$
such that they intersect in three disjoint (singular) conics on the
cubic surface.
We then show that these reducible surfaces can be constructed
directly, in a slightly more general form. In fact we construct a
cubic Del Pezzo surface~$X_0$ and
a smooth general type surface~$T$
of degree~$10$ and sectional genus~$10$
such that~$T \cap X_0$ is the disjoint union of
three smooth conics
and show that the general surface linked $(5,5)$
to~$T \cup X_0$ is
a smooth irregular elliptic surface.
The family of irregular elliptic surfaces constructed via liaison
includes the family of irregular elliptic surfaces constructed via monad
as a special case.

\vspace{3mm}
\noindent
{\it Acknowledgment.}
We thank David Eisenbud, Lucian Badescu, Sorin Popescu and Alfio Ragusa
for organizing the Nato advanced workshop in Erice, September 2001.
Our collaboration started in the stimulating atmosphere of that
meeting.
The first author would like to thank Wolfram Decker 
for many stimulating conversations. 
The first author is also grateful to the DAAD for its financial support.

\section{Preliminaries}
In this part
we recall basic results on
smooth surfaces in~$\P^4$ needed for
the two methods of construction that we shall use.
\subsection{Notation and the basic results}
If not otherwise mentioned,
$X$ denotes a smooth surface in~$\P^4$ and
\begin{itemize}
\item[-] $H = H_X$ its hyperplane class;
\item[-] $K = K_X$ its canonical divisor;
\item[-] $d = d(X) = H^2$ its degree,
\item[-] $\pi = \pi (X)$ its sectional genus (the genus of a general
hyperplane section);
\item[-] $\chi = \chi (\OO_X) = p_g - q + 1$
the Euler-Poincar\'e characteristic of its structure sheaf.
\end{itemize}
The numerical invariants of~$X$ satisfy the {\it double point
formula}
\[
    d^2 - 10d - 5H \cdot K - 2K^2 + 12 \chi = 0
\]
(see, for example~\cite{hartshorne}).
For a curve on~$X$ the arithmetic genus
can be computed by the {\it adjunction formula}
\[
    2p_a (C) - 2 = C^2 + C \cdot K.
\]
The double point formula and the adjunction formula
express the self-intersection number~$K^2$
in terms of~$d$, $\pi$ and~$\chi$:
\begin{eqnarray}
    K^2 = \frac{d^2 - 5d - 10\pi + 12\chi + 10}{2}.
\label{eq:self_intersect}
\end{eqnarray}
For curves~$C$ and~$D$ on~$X$ the adjunction formula yields
the following addition rule for the arithmetic genus:
\begin{eqnarray}
    p_a (C \cup D) = p_a(C) + p_a(D) + C \cdot D -1.
\label{eq:add_genus}
\end{eqnarray}
\subsection{Monads}
\label{monad}
To construct a smooth surface~$X$ in~$\P^4$,
we construct its ideal sheaf~$\II_X$.  Our construction of~$\II_X$
(more generally
coherent sheaves on~$\P^n$) uses monads.
Indeed, monads represent a given coherent sheaf in terms of
simpler vector bundles,
such as direct sum of line bundles,
and homomorphisms between these simpler bundles:
\begin{df}
A {\it monad} for a coherent sheaf~$\GG$ on~$\P^n$
is a complex
\begin{eqnarray}
    \AA \stackrel{\alpha}{\longrightarrow} \BB
    \stackrel{\beta}{\longrightarrow} \CC
\label{eq:monad}
\end{eqnarray}
of vector bundles on~$\P^n$ with~$\alpha$ an injective
map and~$\beta$ a surjective map such that
the cohomology $\ker (\beta) / \im(\alpha)$ is the coherent
sheaf~$\GG$. The {\it display} of the monad~(\ref{eq:monad})
is the following commutative diagram:
\begin{eqnarray}
\begin{CD}
       &  &               &  & 0            & &            0          \\
       &  &               &  & @VVV               @VVV       \\
     0 @>>> \AA @>>> \KK @>>> \GG @>>> 0 \\
       &  &     \|        &  & @VVV               @VVV       \\
     0 @>>> \AA @>{\alpha}>> \BB @>>> \QQ @>>> 0,  \\
       &  &               &  & @VV{\beta}V               @VVV       \\
       &  &               &  & \CC @=  \CC           \\
       &  &               &  & @VVV               @VVV       \\
       &  &               &  & 0           & &    0
\end{CD}
\label{eq:display_1}
\end{eqnarray}
where~$\KK = \ker (\beta)$ and~$\QQ = \coker (\alpha)$.
\end{df}

To find a monad for our ideal sheaf we use Horrocks' technique of
killing cohomology.
First we recall the definition of Horrocks monad:
\begin{df} A monad
\[
    \AA \stackrel{\alpha}{\longrightarrow} \BB
    \stackrel{\beta}{\longrightarrow} \CC
\]
for a coherent sheaf~$\GG$ on~$\P^4$ is a {\it Horrocks monad} if
\begin{enumerate}
    \item[(a)] $\AA$ and $\CC$ are direct sums of line bundles;
    \item[(b)] $\BB$ satisfies
\[
    \H^i_* \BB \simeq
     \begin{cases}
       0 & i=1,n-1 \\
       \H^i_* \GG  & 1 < i < n-1.
     \end{cases}
\]
\end{enumerate}
    We say that this monad is {\it minimal} if
\begin{enumerate}
    \item[(c)] $\CC$ has no direct summand which is the image
    of a line subbundle of $\BB$;
    \item[(d)] there is no direct summand of $\AA$ which is
    isomorphic to a direct summand of $\BB$.
\end{enumerate}
\end{df}

For any vector bundle on~$\P^n$, $n \geq 2$,
there is a Horrocks minimal monad (see, for example~\cite{bh}).
We include a proof of the following special case
for lack of a suitable reference.
\begin{prop}
There is a Horrocks minimal monad for the ideal sheaf~$\II_X$
of any $2$-codimensional locally Cohen-Macaulay
closed subscheme~$X$ of~$\P^n$.
\label{th:horrocks_monad}
\end{prop}

\begin{proof}
Every minimal set of generators for $\H^1_* \II_X$
gives an epimorphism
\begin{eqnarray}
    L_0 \rightarrow \H^1_* \II_X \rightarrow 0,
\label{eq:epimorphism1}
\end{eqnarray}
where~$L_0$ is a free $S = \C [x_0, \dots ,x_n]$-module.
This epimorphism defines an extension
\begin{eqnarray}
    0 \rightarrow \II_X \rightarrow \QQ
    \rightarrow \widetilde{L_0} \rightarrow 0
\label{eq:ex_seq_1}
\end{eqnarray}
of~$\II_X$ by~$\widetilde{L_0}$,
because the given epimorphism~(\ref{eq:epimorphism1}) can
be regarded as an element
of~$\H^1 (\II_X \otimes \widetilde{L_0}) \simeq \Ext^1 (\widetilde{L_0},
\II_X)$.
Taking cohomology~(\ref{eq:ex_seq_1}) gives~$\H^1_* \QQ = 0$
because of~(\ref{eq:epimorphism1}).
Given a minimal set of generators for the
$S$-module~$\H^0_* \omega_X (n+1)$,
we have an epimorphism
\[
    P_0 \rightarrow \H^0_* \omega_X (n+1) \rightarrow 0,
\]
where $P_0$ is a free $S$-module.
Note that~$\omega_X (n+1) = \EE xt^2 (\OO_X, \OO) \simeq
\EE xt^1 (\II_X,\OO)$.
So~$\H^0 (\omega_X (n+1) \otimes \widetilde{P_0^{\vee}})
\simeq \H^0 (\EE xt^1 (\II_X,\OO) \otimes \widetilde{P_0^{\vee}})
\simeq \Ext^1 (\II_X,\widetilde{P_0^{\vee}})$,
and hence we obtain an extension
\begin{eqnarray}
    0 \rightarrow \widetilde{P_0^{\vee}} \rightarrow \KK
    \rightarrow \II_X \rightarrow 0.
\label{eq:ex_seq_2}
\end{eqnarray}
By construction it is clear that~$\H^1_* \KK^{\vee} = 0$,
and hence~$\H^{n-1}_* \KK \simeq \Ext^{n+1}_S (\H^1_* \KK^{\vee},S) = 0$.
The exact sequence~(\ref{eq:ex_seq_1}) induces
\[
     \cdots \rightarrow
    \Ext^1 (\widetilde{L_0}, \widetilde{P_0^{\vee}})
    \rightarrow
    \Ext^1 (\QQ, \widetilde{P_0^{\vee}})
    \rightarrow
    \Ext^1 (\II_X, \widetilde{P_0^{\vee}})
    \rightarrow
    \Ext^2 (\widetilde{L_0}, \widetilde{P_0^{\vee}})
    \rightarrow \cdots.
\]
    Since
    $\Ext^i (\widetilde{L_0}, \widetilde{P_0^{\vee}}) = 0$
    for $i \geq 1$,
    we have an isomorphism
\[
    \Ext^1 (\QQ, \widetilde{P_0^{\vee}})
    \cong
    \Ext^1 (\II_X, \widetilde{P_0^{\vee}}).
\]
    Let us denote by
\[
    0 \rightarrow \widetilde{P_0^{\vee}} \rightarrow \BB
    \rightarrow \QQ \rightarrow 0
\]
    the extension corresponding to~(\ref{eq:ex_seq_2}).
    This extension can be used to complete
    the exact sequences~(\ref{eq:ex_seq_1}) and~(\ref{eq:ex_seq_2})
    to the display
\begin{eqnarray}
    \begin{CD}
       &  &               &  & 0            & &            0          \\
       &  &               &  & @VVV               @VVV       \\
     0 @>>> \widetilde{P_0^{\vee }} @>>> \KK @>>> \II_X @>>> 0 \\
       &  &     \|        &  & @VVV               @VVV       \\
     0 @>>> \widetilde{P_0^{\vee }} @>{\alpha}>> \BB @>>> \QQ @>>> 0 \\
       &  &               &  & @VV{\beta}V               @VVV       \\
       &  &               &  & \widetilde{L_0} @=  \widetilde{L_0}           \\
       &  &               &  & @VVV               @VVV       \\
       &  &               &  & 0           & &    0
\end{CD}
\label{eq:display_2}
\end{eqnarray}
    of a monad
\begin{eqnarray}
    \widetilde{P_0^{\vee}} \stackrel{\alpha}{\longrightarrow} \BB
    \stackrel{\beta}{\longrightarrow}  \widetilde{L_0}
\label{eq:monad_2}
\end{eqnarray}
    for $\II_X$. From the display~(\ref{eq:display_2})
    it follows that $\BB$ satisfies (b),
    because both the modules~$\H^1_* \QQ$ and~$\H^{n-1}_* \KK$
    vanish. So the monad~(\ref{eq:monad_2}) is
    a Horrocks monad.

    Finally we show the minimality of~(\ref{eq:monad_2}).
    It is obvious that
    $\widetilde{L_0}$ (resp. $\widetilde{P_0}$)
    has a direct summand which is the image
    of a line subbundle of $\BB$ (resp. $\BB^{\vee}$) if and only if
    $\QQ$ (resp. $\KK$) has a direct summand which is the image
    of a line subbundle of $\BB$ (resp. $\BB^{\vee}$).
    By construction, however,
    $\QQ$ and $\KK$ do not have such a direct summand,
    which completes the proof.
\end{proof}

\begin{rem}
For the ideal sheaf~$\II_X$ of a smooth surface~$X$ in~$\P^4$
the Horrocks monad is closely related to the minimal free
resolutions of the Hartshorne-Rao modules.
Suppose that~$\H^2_* \II_X$ does not vanish.
Let us denote the minimal free resolution of~$\H^2_* \II_X$ by
\[
    0 \longleftarrow \H^2_* \II_X \longleftarrow
    F_0 \stackrel{f_0}{\longleftarrow}
    F_1 \stackrel{f_1}{\longleftarrow}
    F_2 \stackrel{f_2}{\longleftarrow}
    F_3 \stackrel{f_3}{\longleftarrow}
    F_4 \stackrel{f_4}{\longleftarrow}
    F_5 \longleftarrow 0.
\]
Then $\BB$ has the same intermediate cohomology modules as
the second sheafified syzygy module of~$\H^2_* \II_X$.
Hence~$\H^0_* \BB$ is a direct sum of~$B = \ker (f_1)$ and
a free $S$-module (see~\cite{des} for the proof).
Frequently, $\BB$ is isomorphic to $\widetilde{B}$,
and the minimal free resolution of~$\H^2_* \II_X$ decomposes as
$$
\xymatrix{
\cdots & \ar[l] F_1 & \ar[l] L_0 \oplus F_2'
&&\ar[ll]_{\left(
        \begin{array}{cc}
          \alpha_0 & * \\
          0 & *
        \end{array}
       \right)}
\ar[dl] L_1 \oplus F_3'
& \ar[l] F_4 &\ar[l]\cdots \\
&&& \ar[ul] \ar[dl] B  \\
&&0  && \ar[ul] 0
}
$$

\noindent
including the minimal free presentation of~$\H^1_* \II_X$:
\[
    0 \longleftarrow \H^1_* \II_X \longleftarrow L_0
    \stackrel{\alpha_0}{\longleftarrow} L_1.
\]
So to construct a smooth surface~$X$ in~$\P^4$,
we construct the Hartshorne-Rao modules~$\H^1_* \II_X$ and
$\H^2_* \II_X$ first.

In this case the bundle~$\KK$ can be described explicitly.
So to construct the ideal sheaf~$\II_X$, we use
the determinantal construction following~\cite{des},
namely we establish that the determinantal locus
of the general morphism~$\psi \in \Hom (\widetilde{P_0^{\vee}},\KK)$
is indeed a smooth surface.
\end{rem}

To construct~$\H^1_* \II_X$ and~$\H^2_* \II_X$,
we need information on the dimensions~$\h^i \II_X(m)$
in some range of twists.
This information can be obtained from the Riemann-Roch theorem:

\begin{thm}[Riemann-Roch]
    Let $X$ be a smooth surface in $\P^4$ of degree $d$,
    sectional genus $\pi$, geometric genus $p_g$ and
    irregularity $q$.
    Then
\[
    \chi (\II_X(m)) = \chi (\OO(m)) -
    \left(\begin{matrix} m+1 \\ 2 \end{matrix} \right)d + m(\pi -1)
    -1 + q - p_g.
\]
\label{th:riemann_roch}
\end{thm}

Since surfaces on quadric and cubic hypersurfaces are completely
classified (cf.~\cite{koelblen}), we may assume that
\[
    \h^0 \II_X (2) = \h^0 \II_X(3) = 0.
\]
We observe that some of the cohomology groups vanish:

\begin{prop}[\cite{des}]
    Let $X$ be a smooth non-general type surface in $\P^4$
    which is not lying on any cubic hypersurface.
    Then we have the following table for the dimensions $\h^i \II_X(j)$:
\newlength{\wi}
\newlength{\hi}
{
\fontsize{10pt}{11pt}
\selectfont
$$
{
\setlength{\wi}{14mm}
\setlength{\hi}{8mm}
\fontsize{10pt}{11pt}
\selectfont
\begin{xy}
%
%
<-2\wi,0\hi>;<4\wi,0\hi>**@{-}?>*@{>}?(0.9)*!/^3mm/{j}
%
%
,<-1\wi,0\hi>;<-1\wi,6\hi>**@{-}?>*@{>}?(0.9)*!/^3mm/{i}
%
%
,0+<-2\wi,5\hi>;<3\wi,5\hi>**@{-}
,0+<-2\wi,4\hi>;<3\wi,4\hi>**@{-}
,0+<-2\wi,3\hi>;<3\wi,3\hi>**@{-}
,0+<-2\wi,2\hi>;<3\wi,2\hi>**@{-}
,0+<-2\wi,\hi>;<3\wi,\hi>**@{-}
%
%
,0+<-2\wi,0\hi>;<-2\wi,5\hi>**@{-}
,0+<3\wi,0\hi>;<3\wi,5\hi>**@{-}
,0+<2\wi,0\hi>;<2\wi,5\hi>**@{-}
,0+<1\wi,0\hi>;<1\wi,5\hi>**@{-}
,0+<0\wi,0\hi>;<0\wi,5\hi>**@{-}
,0+<-1\wi,0\hi>;<-1\wi,5\hi>**@{-}
%
%
,0+<-1.5\wi,4.5\hi>*{0}
,0+<-1.5\wi,3.5\hi>*{N+1}
,0+<-1.5\wi,2.5\hi>*{0}
,0+<-1.5\wi,1.5\hi>*{0}
,0+<-1.5\wi,0.5\hi>*{0}
,0+<-.5\wi,4.5\hi>*{0}
,0+<-.5\wi,3.5\hi>*{p_g}
,0+<-.5\wi,2.5\hi>*{q}
,0+<-.5\wi,1.5\hi>*{0}
,0+<-.5\wi,0.5\hi>*{0}
,0+<.5\wi,4.5\hi>*{0}
,0+<.5\wi,3.5\hi>*{0}
,0+<.5\wi,2.5\hi>*{s}
,0+<.5\wi,1.5\hi>*{0}
,0+<.5\wi,0.5\hi>*{0}
,0+<1.5\wi,4.5\hi>*{0}
,0+<1.5\wi,3.5\hi>*{0}
,0+<1.5\wi,2.5\hi>*{\h^2 \II_X(2)}
,0+<1.5\wi,1.5\hi>*{\h^1 \II_X(2)}
,0+<1.5\wi,0.5\hi>*{0}
,0+<2.5\wi,4.5\hi>*{0}
,0+<2.5\wi,3.5\hi>*{0}
,0+<2.5\wi,1.5\hi>*{\h^1 \II_X(3)}
,0+<2.5\wi,2.5\hi>*{\h^2 \II_X(3)}
,0+<2.5\wi,.5\hi>*{0}
,0+<4.5\wi,2.5\hi>*{\h^i \II_X(j)}
%
%
,<-1.5\wi,-.5\hi>*{-1}
,<-.5\wi,-.5\hi>*{0}
,<.5\wi,-.5\hi>*{1}
,<1.5\wi,-.5\hi>*{2}
,<2.5\wi,-.5\hi>*{3}
%
%
,<3.2\wi,.5\hi>*{0}
,<3.2\wi,1.5\hi>*{1}
,<3.2\wi,2.5\hi>*{2}
,<3.2\wi,3.5\hi>*{3}
,<3.2\wi,4.5\hi>*{4}
\end{xy}
}
$$
}

\noindent
where $N = \pi - q + p_g -1$ and $s = \pi - d + 3 + q - p_g$.
\label{th:cohomology_table}
\end{prop}

\noindent
In the sequel we represent a zero in a cohomology table
by an empty box.
\subsection{Liaison}
\label{liaison}
We recall the definition and some basic results from~\cite{ps}.
\begin{df}
Let $X$ and $X'$ be surfaces in~$\P^4$
with no irreducible component in common.
Two surfaces $X$ and $X'$ are said to be {\it linked} $(m,n)$
if there exist hypersurfaces~$V$ and~$V'$ of degree~$m$
and~$n$ respectively such that~$V \cap V' = X \cup X'$.
\end{df}

\noindent
There are two standard sequences of linkage, namely
\[
    0 \rightarrow \OO_X (K) \rightarrow \OO_{X \cup X'}(m+n-5)
    \rightarrow \OO_{X'}(m+n-5) \rightarrow 0
\]
and
\[
    0 \rightarrow \OO_X (K) \rightarrow \OO_{X}(m+n-5)
    \rightarrow \OO_{X \cap X'}(m+n-5) \rightarrow 0.
\]

For the construction we will use the following theorem:
\begin{prop}[\cite{ps}]
    Let $X$ be a locally complete intersection surface in~$\P^4$.
    If $X$ is scheme-theoretically cut out by
    hypersurfaces of degree~$d$,
    then $X$ is linked to a smooth surface~$X'$
    in the complete intersection of two hypersurfaces
    of degree~$d$.
\label{th:smoothness}
\end{prop}

\section{A monad construction}
Let $\P^4$ be the $4$-dimensional projective space over $\C$
with homogeneous coordinate ring~$S = \C [x_0, \dots ,x_4]$.
We construct an example of a smooth irregular
proper elliptic surface in~$\P^4$ using Horrocks' technique of
killing cohomology.
\begin{thm}
    There exists a smooth, minimal, proper elliptic surface in~$\P^4$
    with $d =~12$, $\pi =~13$, $p_g =~3$ and~$q = 1$.
\label{th:monad_const}
\end{thm}
\begin{proof}
    For a smooth surface $X$ in~$\P^4$ with
    the given invariants, formula~(\ref{eq:self_intersect})
    leads to~$K^2 = 0$.
    By the Riemann-Roch formula and Proposition~\ref{th:cohomology_table},
    the table
\newlength{\br}
\newlength{\ho}
{
\fontsize{10pt}{11pt}
\selectfont
$$
{
\setlength{\br}{10mm}
\setlength{\ho}{6mm}
\fontsize{10pt}{11pt}
\selectfont
\begin{xy}
%
%
<-2\br,0\ho>;<4\br,0\ho>**@{-}?>*@{>}?(0.9)*!/^3mm/{j}
%
%
,<-1\br,0\ho>;<-1\br,6\ho>**@{-}?>*@{>}?(0.9)*!/^3mm/{i}
%
%
,0+<-2\br,5\ho>;<3\br,5\ho>**@{-}
,0+<-2\br,4\ho>;<3\br,4\ho>**@{-}
,0+<-2\br,3\ho>;<3\br,3\ho>**@{-}
,0+<-2\br,2\ho>;<3\br,2\ho>**@{-}
,0+<-2\br,\ho>;<3\br,\ho>**@{-}
%
%
,0+<-2\br,0\ho>;<-2\br,5\ho>**@{-}
,0+<3\br,0\ho>;<3\br,5\ho>**@{-}
,0+<2\br,0\ho>;<2\br,5\ho>**@{-}
,0+<1\br,0\ho>;<1\br,5\ho>**@{-}
,0+<0\br,0\ho>;<0\br,5\ho>**@{-}
,0+<-1\br,0\ho>;<-1\br,5\ho>**@{-}
%
%
,0+<-1.5\br,3.5\ho>*{15}
,0+<-.5\br,3.5\ho>*{3}
,0+<-.5\br,2.5\ho>*{1}
,0+<.5\br,2.5\ho>*{2}
,0+<1.5\br,2.5\ho>*{a}
,0+<1.5\br,1.5\ho>*{a}
,0+<2.5\br,1.5\ho>*{b+4}
,0+<2.5\br,2.5\ho>*{b}
,0+<4.5\br,2.5\ho>*{\h^i \II_X(j)}
%
%
,<-1.5\br,-.5\ho>*{-1}
,<-.5\br,-.5\ho>*{0}
,<.5\br,-.5\ho>*{1}
,<1.5\br,-.5\ho>*{2}
,<2.5\br,-.5\ho>*{3}
%
%
,<3.2\br,.5\ho>*{0}
,<3.2\br,1.5\ho>*{1}
,<3.2\br,2.5\ho>*{2}
,<3.2\br,3.5\ho>*{3}
,<3.2\br,4.5\ho>*{4}
\end{xy}
}
$$
}

\noindent
    reflects the relations between the dimensions~$\h^i \II_X(j)$
    in the range~$-1 \leq j \leq 3$ of twists.
    A simple example of a finite length graded $S$-module with
    Hilbert function $(1,2,\cdots )$ is
    the module~$M$ with minimal free presentation:
\[
    0 \longleftarrow M \longleftarrow S \stackrel{f_0}{\longleftarrow}
    3S(-1) \oplus 2S(-2),
\]
    where $f_0 = (x_0,x_1,x_2,x_3^2,x_4^2)$.
    The minimal free resolution of~$M$ is of the following form:
\[
\vbox{%
\halign{&\hfil\,$#$\,\hfil\cr
0&\leftarrow&M&\leftarrow&S&&3S(-1)&&3S(-2)&&S(-3)\cr
&&&&&\vbox to 10pt{\vskip-4pt\hbox{$\nwarrow$}\vss}&\oplus&&\oplus&&\oplus \cr
&&&&&&2S(-2)&\vbox to 10pt{\vskip-4pt\hbox{$\nwarrow$}\vss}
&6S(-3)&\leftarrow&6S(-4)&&2S(-5)\cr
&&&&&&&&\oplus&&\oplus&\vbox to 10pt{\vskip-4pt\hbox{$\nwarrow$}\vss}
&\oplus& \cr
&&&&&&&&S(-4)&&3S(-5)&&3S(-6)&\vbox to 10pt{\vskip-4pt\hbox{$\nwarrow$}\vss}
&S(-7)&\leftarrow&0.  \cr
}}
\]
    Suppose that $\H^2_* \II_X = M$.
    Then~$a = 1$ and~$b = 0$.
    Let~$B$ be the second syzygy module of~$M$.
    We will suppose that~$\widetilde{B}$ is the middle
    term of the Horrocks minimal monad for~$\II_X$:
\[
    \widetilde{P_0^{\vee}} \rightarrow \widetilde{B}
    \rightarrow \widetilde{L_0}.
\]
    With suitably chosen bases,
    $B$ is generated by the column of the matrix
\[
    f_2 =
    \left(
     \begin{matrix}
     x_2 &
     {{x}}_{{3}}^{2}&
        0&
        0&
        {{x}}_{{4}}^{2}&
        0&
        0&
        0&
        0&
        0\\
        {-{{x}}_{1}}&
        0&
        {{x}}_{{3}}^{2}&
        0&
        0&
        {{x}}_{{4}}^{2}&
        0&
        0&
        0&
        0\\
        {{x}}_{0}&
        0&
        0&
        {{x}}_{{3}}^{2}&
        0&
        0&
        {{x}}_{{4}}^{2}&
        0&
        0&
        0\\
        0&
        {-{{x}}_{1}}&
        {-{{x}}_{{2}}}&
        0&
        0&
        0&
        0&
        {{x}}_{{4}}^{2}&
        0&
        0\\
        0&
        {{x}}_{0}&
        0&
        {-{{x}}_{{2}}}&
        0&
        0&
        0&
        0&
        {{x}}_{{4}}^{2}&
        0\\
        0&
        0&
         {{x}}_{0}&
        {{x}}_{1}&
        0&
        0&
        0&
        0&
        0&
        {{x}}_{{4}}^{2}\\
        0&
        0&
        0&
        0&
        {-{{x}}_{1}}&
        {-{{x}}_{{2}}}&
        0&
        {-{{x}}_{{3}}^{2}}&
        0&
        0\\
        0&
        0&
        0&
        0&
        {{x}}_{0}&
        0&
        {-{{x}}_{{2}}}&
        0&
        {-{{x}}_{{3}}^{2}}&
        0\\
        0&
        0&
        0&
        0&
        0&
        {{x}}_{0}&
        {{x}}_{1}&
        0&
        0&
        {-{{x}}_{{3}}^{2}}\\
        0&
        0&
        0&
        0&
        0&
        0&
        0&
        {{x}}_{0}&
        {{x}}_{1}&
        {{x}}_{{2}}\\
      \end{matrix}
    \right).
\]
    For simplicity we denote by~$F_2$ and~$F_3$
    the target and source of~$f_2$ respectively.
    From the minimal free resolution of~$M$ it follows that
    $\widetilde{B}$ has the following
    Beilinson cohomology table:
{
\fontsize{10pt}{11pt}
\selectfont
$$
{
\setlength{\br}{10mm}
\setlength{\ho}{6mm}
\fontsize{10pt}{11pt}
\selectfont
\begin{xy}
%
%
<-2\br,0\ho>;<4\br,0\ho>**@{-}?>*@{>}?(0.9)*!/^3mm/{j}
%
%
,<-1\br,0\ho>;<-1\br,6\ho>**@{-}?>*@{>}?(0.9)*!/^3mm/{i}
%
%
,0+<-2\br,5\ho>;<3\br,5\ho>**@{-}
,0+<-2\br,4\ho>;<3\br,4\ho>**@{-}
,0+<-2\br,3\ho>;<3\br,3\ho>**@{-}
,0+<-2\br,2\ho>;<3\br,2\ho>**@{-}
,0+<-2\br,\ho>;<3\br,\ho>**@{-}
%
%
,0+<-2\br,0\ho>;<-2\br,5\ho>**@{-}
,0+<3\br,0\ho>;<3\br,5\ho>**@{-}
,0+<2\br,0\ho>;<2\br,5\ho>**@{-}
,0+<1\br,0\ho>;<1\br,5\ho>**@{-}
,0+<0\br,0\ho>;<0\br,5\ho>**@{-}
,0+<-1\br,0\ho>;<-1\br,5\ho>**@{-}
%
%
,0+<-1.5\br,4.5\ho>*{1}
,0+<-1.5\br,2.5\ho>*{}
,0+<-.5\br,2.5\ho>*{1}
,0+<-.5\br,2.5\ho>*{}
,0+<.5\br,2.5\ho>*{2}
,0+<.5\br,1.5\ho>*{}
,0+<1.5\br,2.5\ho>*{1}
,0+<1.5\br,1.5\ho>*{}
,0+<2.5\br,1.5\ho>*{}
,0+<2.5\br,.5\ho>*{1}
,0+<4.5\br,2.5\ho>*{\h^i \widetilde{B}(j)}
%
%
,<-1.5\br,-.5\ho>*{-1}
,<-.5\br,-.5\ho>*{0}
,<.5\br,-.5\ho>*{1}
,<1.5\br,-.5\ho>*{2}
,<2.5\br,-.5\ho>*{3}
%
%
,<3.2\br,.5\ho>*{0}
,<3.2\br,1.5\ho>*{1}
,<3.2\br,2.5\ho>*{2}
,<3.2\br,3.5\ho>*{3}
,<3.2\br,4.5\ho>*{4}
\end{xy}
}
$$
}

\noindent
    Suppose that~$\H^1_* \II_X$ is generated
    in the first non-zero twist, that is, monogeneous.
    Then~$\widetilde{L_0}~=~\OO(-2)$.
    Since~$p_g = \h^3 \II_X = 3$,
    the map~$3S \rightarrow \H^0_* \omega_X$ forms
    part of the minimal generating set of~$\H^0_* \omega_X$.
    From the cohomology table of~$\widetilde{B}$,
\[
   \h^4 \widetilde{P_0^{\vee}}(-1) = \h^4 \II_X(-1) + h^4 \widetilde{B}(-1)
   = 15 + 1 = 16.
\]
    So~$\widetilde{P_0^{\vee}}$ must contain~$\OO(-4)$
    as a direct summand.
    Since~$\rank (\widetilde{B}) = 6$,
    we can deduce that~$\widetilde{P_0^{\vee}} \simeq 3\OO(-5) \oplus \OO(-4)$.
    So the Horrocks monad for~$\II_X$ is of type
\[
    3\OO(-5) \oplus \OO(-4) \rightarrow \widetilde{B}
    \rightarrow \OO(-2)
\]
    and we have the display of the monad:
\begin{eqnarray}
    \begin{CD}
       &  &               &  & 0            & &            0          \\
       &  &               &  & @VVV               @VVV       \\
     0 @>>> 3\OO(-5) \oplus \OO(-4) @>>> \KK @>>> \II_X @>>> 0 \\
       &  &     \|        &  & @VVV               @VVV       \\
     0 @>>> 3\OO(-5) \oplus \OO(-4) @>>> \BB @>>> \QQ @>>> 0.  \\
       &  &               &  & @VVV               @VVV       \\
       &  &               &  & \OO(-2) @=  \OO(-2)    \\
       &  &               &  & @VVV               @VVV       \\
       &  &               &  & 0           & &    0
\end{CD}
\label{eq:display_3.1}
\end{eqnarray}

    Next we look for an appropriate finite length
    graded module~$N$ with Hilbert function~$(1,4, \cdots)$
    which is monogeneous.
    Recall that $P_0^{\vee}$ contains~$S(-4)$ as a direct summand.
    So to construct such a module,
    we have to find a map~$\phi : F_2 \rightarrow S(-2)$
    such that the source of the syzygy matrix
    of~$f_2 \circ \phi : F_3 \rightarrow S(-2)$
    includes~$S(-4)$ as a direct summand.
    For instance,
    the map given by the matrix
\[
    \phi =
    \left(
       1 , 0 , 0 , 0 , 0 , -x_0 , 0 , 0 , x_1 , 0
     \right) : F_2 \rightarrow S(-2)
\]
    satisfies the desired condition.
    One shows by straightforward calculations that
    the cokernel of~$f_2 \circ \phi$
    is of finite length and that~$N$ has the minimal free resolution
    of type
\[
\vbox{%
\halign{&\hfil\,$#$\,\hfil\cr
0&\leftarrow&N&\leftarrow&S(-2)&&S(-3)\cr
&&&&&\vbox to 8pt{\vskip-3pt\hbox{$\nwarrow$}\vss}&\oplus \cr
&&&&&&5S(-4)&\vbox to 8pt{\vskip-3pt\hbox{$\nwarrow$}\vss}&7S(-5)&&2S(-6) \cr
&&&&&&&&\oplus&&\oplus \cr
&&&&&&&&7S(-6)&\vbox to 8pt{\vskip-3pt\hbox{$\nwarrow$}\vss}&11S(-7)&&4S(-8)\cr
&&&&&&&&&&\oplus&\vbox to 8pt{\vskip-3pt\hbox{$\nwarrow$}\vss}&\oplus\cr
&&&&&&&&&&3S(-8)&&5S(-9)&\vbox to
8pt{\vskip-3pt\hbox{$\nwarrow$}\vss}&2S(-10)&\leftarrow&0. \cr
}}
\]

\noindent
In this case~$\KK$ has syzygies of type
\[
\vbox{%
\halign{&\hfil\,$#$\,\hfil\cr
&&&&\OO(-4)\cr
&&&&\oplus\cr
0&\leftarrow&\KK&\leftarrow&8\OO(-5)&&&2\OO(-6)\cr
&&&&\oplus&&& \oplus \cr
&&&& 4\OO (-6)&&\vbox to 10pt{\vskip-4pt\hbox{$\nwarrow$}\vss}
&10\OO(-7)&&4 \OO(-8)\cr
&&&&&&&\oplus  &\vbox to 10pt{\vskip-4pt\hbox{$\nwarrow$}\vss}& \oplus \cr
&&&&&&&3 \OO(-8) && 5\OO(-9)
&\vbox to 10pt{\vskip-4pt\hbox{$\nwarrow$}\vss}& 2\OO (-10)&\leftarrow&0\ .\cr
}}
\]

\noindent
The ideal sheaf~$\II_X$ will be obtained via an exact sequence
\[
    0 \rightarrow 3\OO(-5) \oplus \OO(-4) \rightarrow \KK
    \rightarrow \II_X \rightarrow 0.
\]
    However~$X$ will be obtained as the determinantal locus of
    the general map between the vector bundles $3\OO(-5)$ and~$\EE$,
    where~$\EE$ is the cokernel of
    a non-zero map~$\OO(-4) \rightarrow \KK$.
    Indeed, the map
\[
    \psi =
    \mbox{$^t\left(
     0,0,0,1,0,1,0,0,0,0 \right)$} : S(-4) \rightarrow F_3
\]
    forms part of the minimal generating set of~$\Syz (f_2 \circ \phi)$
    and a straightforward calculation shows that
    the cokernel of~$(\psi \circ f_2)^{\vee}$
    is of finite length.
    This means that for a non-zero element of~$\Hom (\OO(-4), \KK)$
    the cokernel~$\EE$ is a rank four vector bundle.
    Clearly the minimal free resolution of~$\EE$ is of the form
\[
\vbox{%
\halign{&\hfil\,$#$\,\hfil\cr
&&&&8\OO(-5)&&&2\OO(-6)\cr
0&\leftarrow&\EE&\leftarrow&\oplus&&& \oplus \cr
&&&& 4\OO (-6)&&\vbox to 10pt{\vskip-4pt\hbox{$\nwarrow$}\vss}
&10\OO(-7)&&4 \OO(-8)\cr
&&&&&&&\oplus  &\vbox to 10pt{\vskip-4pt\hbox{$\nwarrow$}\vss}& \oplus \cr
&&&&&&&3 \OO(-8) && 5\OO(-9)
&\vbox to 10pt{\vskip-4pt\hbox{$\nwarrow$}\vss}& 2\OO (-10)&\leftarrow&0\ .\cr
}}
\]
 By using {\it Macaulay $2$}, we can check that
 the determinantal locus~$V(f)$ of the general
 map~$f\in\Hom~(3\OO(-5),\EE)$
 is a smooth surface. 
 The invariants are computed from the minimal free resolution 
 of~$\II_X$. 
 Since~$\II_X$ has the presentation  
\begin{eqnarray}
    0 \rightarrow 3\OO(-5) \rightarrow \EE \rightarrow \II_X
    \rightarrow 0, 
\label{eq:exact_seq_3.1}
\end{eqnarray}
 the minimal free resolution of $\II_X$ can be computed
 from that of~$\EE$ and has the shape 
\[
\vbox{%
\halign{&\hfil\,$#$\,\hfil\cr
&&&&5\OO(-5)&&&2\OO(-6)\cr
0&\leftarrow&\II_X&\leftarrow&\oplus&&& \oplus \cr
&&&& 4\OO (-6)&&\vbox to 10pt{\vskip-4pt\hbox{$\nwarrow$}\vss}
&10\OO(-7)&&4 \OO(-8)\cr
&&&&&&&\oplus  &\vbox to 10pt{\vskip-4pt\hbox{$\nwarrow$}\vss}& \oplus \cr
&&&&&&&3 \OO(-8) && 5\OO(-9)
&\vbox to 10pt{\vskip-4pt\hbox{$\nwarrow$}\vss}& 2\OO (-10)&\leftarrow&0\ .\cr
}}
\]
 Computing the Hilbert Polynomial of~$X$ from this resolution, 
 we conclude that~$X$ has degree~$d = 12$, sectional genus~$\pi = 13$, 
 and Euler characteristic~$\chi = 3$. From~(\ref{eq:self_intersect}) 
 we now get~$K^2 = 0$. 
 Dualizing~(\ref{eq:exact_seq_3.1}), we obtain
\[
    0 \rightarrow \OO(-5) \rightarrow \EE^{\vee} \rightarrow
    3\OO \rightarrow \omega_X \rightarrow 0.
\]
 So~$p_g = 3$ and $X$ is irregular with~$q = 1$. 
 Furthermore $\omega_X$ is globally generated by its own three sections,
 and thus these sections define the canonical map
\[
    \Phi = \Phi_{\mid K \mid} : X \rightarrow \P^2.
\]
 Then we have the Stein factorization of $\Phi$:
$$
\xymatrix{
X
\ar[dr]_p
\ar[rr]^{\Phi}
&& C \subset \P^2,  \\
& B \ar[ur] \\
}
$$
    where $p$ has connected fibers.
    Let $F$ be a fiber of $p$.
    Then $F \cdot (F+K) = 2K^2 = 0$.
    So it follows from the adjunction formula that~$p_a (F) = 1$,
    which implies that
    the smooth fibers of~$p$ are elliptic curves.
    Therefore~$X$ is an elliptic surface. Finally, since $K^2 = 0$, the
    surface is minimal,
    which completes the proof.
\end{proof}

\begin{rem}
    The surface constructed above may also be obtained as 
    the dependency locus of two general sections of a rank
    three vector bundle on~$\P^4$. 
    In fact, in the notation of the above proof, 
    define a homomorphism by
\[
  \psi' = \mbox{$^t (0,x_1,0,0,x_0,0,0,0,0,1)$} : S(-5) \rightarrow F_3. 
\]
    Then $\psi'$ forms part of the minimal set of generators for
    $\Syz (f_2 \circ \phi)$.
    So~$\psi'$ induces a map~$g : \OO(-5) \rightarrow \EE$.
    One can show by a straightforward calculation that
    $g$ is an injective bundle map.
    The cokernel~$\FF$ of~$g$ is a rank three vector bundle on~$\P^4$,
    and thus $X$ can be regarded as
    the dependency locus of two general sections of~$\FF$.
\label{th:rank_three}
\end{rem}

\section{A liaison construction}
    We first motivate this second construction. 
    Let $X$ be the irregular elliptic surface in $\P^4$
    constructed in the previous section.
    With {\it Macaulay~$2$} we can check that
    the quintics containing $X$ intersect in
\[
    V \left(\H^0(\II_X(5))\right) = X \cup X_0,
\]
    where $X_0$ is a singular rational cubic scroll surface
    in a hyperplane $H_0$ of $\P^4$.  The scroll $X_{0}$
    may be obtained by projecting
    a smooth rational cubic scroll
    $S$ in $\P^4$ from a point off $S$. The scroll $S$ contains a unique
    directrix line. Let $L$ be the image of this line on
    $X_{0}$.

    We can check a
    few more facts with {\it Macaulay~$2$}.
    The surface $X$ does not intersect the line $L$.
    On the other hand the hyperplane section $H_0 \cap X$,
    denoted by $C_0$, is equal to~$X_0 \cap X$,
    and hence $\deg (C_0) = 12$ and $\pi (C_0) = 13$.
    Let $\widetilde{C_{0}}$ and $\widetilde{L}$ be the inverse images of
    $C_{0}$ and $L$ respectively on $S$. Then
    $\widetilde{C_{0}} \equiv 3(H_{S} + \widetilde{L})$,
    since $\widetilde{C_{0}}$ does not intersect~$\widetilde{L}$.
    Now,~$C_0 = X \cap X_0=X\cap H_{0}$, so
    $X \cup X_0$ is locally a complete intersection cut out by
    quintics.
    Therefore, by Proposition~\ref{th:smoothness}, $X \cup X_0$ can be 
    linked~$(5,5)$ to a smooth general type
    surface $T$ of degree $d = 10$,
    $\pi = 10$ and $\chi = 4$.
    This implies that $X$ is reobtained as the residual surface
    of~$X_0 \cup T$ in the complete intersection of two quintics
    in~$\H^0 (\II_{X_0 \cup T}(5))$.
    Let $C = X_0 \cap (X \cup T)$.
    Compare the second standard exact sequence
\[
    0 \rightarrow \OO_{X_0}(K_{X_0}) \rightarrow \OO_{X_0} (5)
    \rightarrow \OO_C(5) \rightarrow 0
\]
    with the exact sequence
\[
    0 \rightarrow \OO_{X_0}(-C) \rightarrow \OO_{X_0}
    \rightarrow \OO_C \rightarrow 0.
\]
    Then we can deduce that $C \equiv 5H_{X_0} - K_{X_0} \equiv 6H_{X_0}$.
    On $S$ we get $\widetilde{C} \equiv 6H_{S}$.
    Let $C_1 := X_0 \cap T$ and let $\widetilde{C_{1}}$ be
    its inverse image on $S$.
    Then $\widetilde{C_{1}}\equiv 3(H_{S} - \widetilde{L})$,
    and hence~$\widetilde{C_{1}}$, is the disjoint union of six members of
    the ruling on $S$.  The singular locus of $X_{0}$ is a line, so $C$
    has six double points along this line.  Three of these appear already
    on $C_{0}$ since $p(\widetilde{C_{0}})=10=p({C_{0}})-3$ and are
    precisely the intersection of $C_{0}$ with the singular locus.  The other
    three points are therefore intersections of pairs of rulings in
    $C_{1}$.  Thus $C_{1}$ consists of three disjoint singular conic
    sections that each lie in a plane that contains the line $L$.

    We recall some facts about smooth general type surfaces~$T$ of
    degree $d = 10$,
    $\pi = 10$ and $\chi = 4$ from~\cite{pr}, Proposition~4.18.  First of all,
    any such surface $T$ can be linked $(4,4)$ to the union of
    a degenerate cubic scroll $U_0$ formed by three planes
    and a Del Pezzo surface~$U_1$ of degree $3$
    such that each plane in~$U_0$ intersect
    the hyperplane of~$U_1$ along the same line $L$
    on~$U_1$.
    The line $L$ is the unique $6$-secant line to $T$.
    Every hyperplane through $L$
    intersect $T$ in a section which is contained in a cubic surface.
    Furthermore, the surface~$T$ is embedded by the linear system
\begin{eqnarray}
    |2K - A_1 - A_2 - A_3|
\label{eq:linear_system}
\end{eqnarray}
    in~$\P^4$, where~$A_1$, $A_2$ and~$A_3$ are pairwise disjoint $(-2)$-curves
    embedded as conics that each intersect $L$ in two points.

    We now want to choose an appropriate surface~$T$ together with a
    cubic surface $X_{0}$ such that their union is linked $(5,5)$ to
    a smooth irregular minimal elliptic surface.
    On such a surface $T$ the three conics in $C_{1}=X_{0}\cap T$ all
    have a common
    secant line which coincides with $L$, so we assume that
    $C_{1}=A_1 \cup A_2 \cup A_3$.

    In particular we want $C_{1}$ to be contained in a
    hyperplane. Let~$D$ be the residual curve~$H_{T} - C_1$.
    Then $C_{1}\cdot D=C_{1}\cdot (H_{T}-C_{1})=6+6=12$ and $p_a (D) =
    1$ by~(\ref{eq:add_genus}).
    This means that~$D$ is an elliptic curve of degree~$4$.
    So to construct an irregular elliptic surface $X$,
    we find a smooth general type surface $T$ in $\P^4$
    with $d = 10$, $\pi = 10$ and $\chi = 4$
    and a hyperplane~$H_0$ of~$\P^4$ such that
\begin{list}{-}{}
    \item $H_0$ contains the $6$-secant $L$ to $T$;
    \item $H_0 \cap T = C_1 \cup D$, where $C_1$ is the disjoint
    union of three smooth $(-2)$-conics and
    $D$ is an elliptic curve of degree~$4$ such that~$D \cap L = \emptyset$
    and~$C_1 \cap D$ consists of~$12$ points.
\end{list}
    Since the hyperplane $H_{0}$ of $C_1 \cup D$ contains $L$, the hyperplane
    section lies in a unique cubic surface, which we denote by $U_{1}$.
    On the other hand the curve $C_1$ lies in the cubic surface
    defined by the three planes of the conics in $C_{1}$.
    The general element of~$\H^0 \II_{C_1,H_0}(3)$ therefore defines
    a smooth cubic surface~$X_0$,
    which does not coincide with~$U_1$.
    We will show that the surface linked $(5,5)$ to the reducible
    surface $T \cup X_0$ is a smooth irregular minimal elliptic surface.
    But first we make the explicit construction.

\vspace{3mm}
\noindent
{\bf Construction:}\ \
Let~$H_0 = V(h_0)$ be a hyperplane of~$\P^4$ and
$L = V(h_0,h_1,h_2)$ a line in~$\P^4$.
Then three quadric minors~$Q_1,Q_2$ and~$Q_3$  of a~$2 \times 3$
matrix whose entries are~$\C$-linear
combinations of~$h_0,h_1$ and~$h_2$, define a degenerate cubic scroll~$U_0$.
By construction, $U_{0}$ then consists of three planes through the line~$L$.
Next, we define a cubic Del Pezzo surface~$U_1$ in~$H_0$
given by a general cubic~$f$ in~$\H^0 \II_{L,H_0}(3)$.
An elliptic curve~$D$ on~$U_1$,
which does not meet~$L$ can be obtained as the general member
of the linear system~$|H_{U_1}+L|$.
Explicitly,
we consider a plane conic~$A \in |H_{U_1}-L|$ on~$U_1$
and take a quadric~$q \in \H^0 \II_{A,H_0}(2)$. Then $D = B - A$,
where~$B = V(h_0,q)$.

The homogeneous ideal~$I_D$ of~$D$ can be written
as
\[
    I_D = (h_0,q_1,q_2),
\]
where~$q_1,q_2 \in \H^0 \OO_{H_0}(2)$ such that
$f \in (q_1,q_2)$.
Consider the double structure~$D_2 = V(h_0^2,q_1,q_2) \subset \P^4$.
Then~$D_2$ is not contained in~$U_1$.
Multiplication by~$h_0$ defines
the following exact sequence:
\[
    0 \longrightarrow  \II_{U_0 \cup D}(3)
    \stackrel{\cdot h_0}{\longrightarrow} \II_{U_1 \cup U_0 \cup D_2}(4)
    \longrightarrow \II_{(U_1 \cup U_0 \cup D_2)\cap H_0,H_0}(4)
    \longrightarrow 0.
\]
Consider the exact sequence
\[
    0 \longrightarrow  \II_{U_0}(2)
    \stackrel{\cdot h_0}{\longrightarrow} \II_{U_0 \cup D}(3)
    \longrightarrow \II_{(U_0 \cup D)\cap H_0,H_0}(3)
    \longrightarrow 0.
\]
The curve $(U_0 \cup D)\cap H_0$ is the union of the
first order neighborhood on~$L$ and the elliptic curve~$D$
of degree~$4$ with~$D \cap L = \emptyset$, so
we have~$\H^0 \II_{(U_0 \cup D)\cap H_0,H_0}(3) = 0$.
We therefore conclude that~$\h^0 \II_{U_0 \cup D}(3) = \h^0 \II_{U_0}(2) = 3$.
Furthermore, $\H^0 \II_{(U_1 \cup U_0 \cup D)\cap H_0,H_0}(4)
\simeq \H^0 \II_{L,H_0}(1)$,
and hence~$\h^0 \II_{U_1 \cup U_0 \cup D_2}(4) \leq 5$.

On the other hand, the following example tells us
that~$\h^0 \II_{U_1 \cup U_0 \cup D_2}(4) \geq 5$:
\begin{ex}
    Let~$S = \C [x_0,\dots ,x_4]$ be the
    homogeneous coordinate ring of~$\P^4$.
    We define the homogeneous ideals~$I_{H_0},I_L,I_{U_0},I_{U_1}$
    and~$I_D$ by
\[
\begin{aligned}
    I_{H_0} &= (h_0) = (x_1); \\
    I_L \   &= (h_0,h_1,h_2) = (x_0,x_1,x_4); \\
    I_{U_0} &= (Q_1,Q_2,Q_3) =
((x_1-x_4)(x_4-x_0),(x_4-x_0)(x_1+x_0),(x_1+x_0)(x_1-x_4)); \\
    I_{U_1} &= (h_0,f) = (x_1,(x_4x_0-x_2x_3)x_0+(x_2^2+x_3^2+x_4^2)x_4); \\
    I_D \   &= (h_0,q_1,q_2) = (x_1,x_4x_0-x_2x_3,x^2+x_3^2+x_4^2).
\end{aligned}
\]
    Then the following quartics form part of the minimal generating
    set of~$I_{U_1 \cup U_0 \cup D_2}$:
\[
\begin{aligned}
    g_1 &= (x_1^2-x_4^2)h_0^2; \\
    g_2 &= (-x_1^2+x_4^2+Q_3)h_0^2; \\
    g_3 &= h_0^2q_2 - (x_0x_4+x_1x_4)Q_1 + (-q_1+x_1x_4-x_4^2)Q_2 +
q_2Q_3 + x_4^2q_2;\\
    g_4 &= (x_1^2-x_4^2 - Q_2 - Q_3)h_0^2; \\
    g_5 &= q_2h_0^2 + (-x_1x_4+x_4^2)Q_2 + (-x_2x_3+x_4^2)Q_3 - x_4^2q_2.
\end{aligned}
\]
\end{ex}

The sequence of global sections
\[
    0 \longrightarrow  \H^0\II_{U_0 \cup D}(3)
    \stackrel{\cdot h_0}{\longrightarrow} \H^0\II_{U_1 \cup U_0 \cup D_2}(4) 
    \longrightarrow \H^0\II_{(U_1 \cup U_0 \cup D_2)\cap H_0,H_0}(4)
    \longrightarrow 0.
\]
is therefore exact.
Let $f_{1},f_{2}$ be general elements in $\H^0 \II_{U_1 \cup U_0 \cup
D_{2}}(4)$.  By Bertini, the complete intersection~$V(f_1)$ and~$V(f_2)$
is singular at most along the base locus of the quartics in $\H^0
\II_{U_1 \cup U_0 \cup D_2}(4)$.  Therefore
$U_0 \cup U_1$ is linked~$(4,4)$ to a smooth surface~$T$
of~$d = 10$ and~$\pi = 10$ in the complete intersection~$V(f_1,f_{2})$.
The surface~$T$ contains~$D_2$, and hence~$D$.
Notice that~$T \cap H_0 = T \cap U_1$.
We denote~$T \cap U_1$ by~$C$
and the residual curve~$C-D$ by~$C_1$.
Since~$\deg (C) = 10$ and~$p_a(C) = 10$,
the divisor~$C$ of~$U_1$ can be written as
\[
    C = 4H_{U_1} - 2L.
\]
So~$C_1 = 3 (H_{U_1}-L)$,
which implies that~$C_1$ is the disjoint union of three conics and the
intersection~$D\cap C_1$ is~$12$ points.
Each conic of~$C_1$ is contained in a unique plane,
and these three planes form a cubic surface~$U_{\infty}$ in~$H_0$.
So we have the pencil of cubics containing~$C_1$
spanned by~$U_1$ and~$U_{\infty}$,
and the general member of this pencil is smooth.
Therefore we  may choose a smooth cubic surface~$X_0$
containing~$C_1$ such that~$X_0 \not= U_1$ and does not contain $D$.
But $D$ is a quartic curve, so $X_0 \cap D $ is 12 points (distinct by
a general choice of $X_{0}$).
Therefore~$X_0 \cap D = C_1 \cap D$,
and~$C_1 = X_0 \cap T$.
In particular~$X_0$ and~$T$ intersect transversally along~$C_1$,
and hence~$Y = T \cup X_0$ is locally complete intersection.
\begin{thm}
The reducible surface~$Y$ is linked~$(5,5)$
to a minimal smooth elliptic surface with~$d = 12$,~$\pi = 13$,~$p_g = 3$, 
and~$q = 1$.
\end{thm}

\begin{proof}
For smoothness it suffices, 
by Proposition~\ref{th:smoothness},
to show that $\II_Y (5)$ is globally generated.
We first estimate~$\h^0 (\II_Y(5))$:
\begin{lem}\label{exact}
    The residual exact sequence
\[
    0 \longrightarrow \II_T(4) \stackrel{\cdot h_0}{\longrightarrow}
    \II_Y (5) \longrightarrow \II_{Y \cap H_0,H_0}(5)
    \longrightarrow 0
\]
    remains exact on global sections.
\end{lem}
\begin{proof}
    We have already shown that the residual curve~$D = C - C_1$ is an
    elliptic curve of degree~$4$.
    So~$\H^0 \II_{Y \cap H_0,H_0}(5) \simeq \H^0 \II_{D\cup X_0, H_0}(5)
    \simeq \H^0 \II_{D,H_0}(2)$,
    and thus~$\h^0 \II_{Y \cap H_0,H_0}(5) = 2$.
    On the other hand~$\h^0 \II_T(4) = 3$
    and~$\h^1 \II_T(4) = 1$ by \cite{pr}.
    In order to prove the lemma,
    it is therefore enough to show that~$\h^0 \II_Y(5) = 5$.
    Let~$H = V(h)$ be a hyperplane containing
    the $6$-secant line~$L$ to~$T$, but not~$X_0$.
    Then multiplication by~$h$ defines
    the exact sequence
\[
     0 \longrightarrow \II_Y(4) \stackrel{\cdot h}{\longrightarrow}
    \II_Y (5) \longrightarrow \II_{Y \cap H,H}(5)
    \longrightarrow 0.
\]
    If~$\h^0 \II_{Y \cap H,H}(5) = 7$, then~$\h^0 \II_Y(5) = 5$,
    because~$\h^0 \II_Y (4) = 0, \h^1 \II_Y(4) = 2$
    and~$\h^2 \II_Y (4) = 0$.
    So we shall show that~$\h^0 \II_{Y \cap H,H}(5) = 7$.

    Since~$H$ contains~$L$,
    the hyperplane section $Y \cap H$ consists of~$L$, a conic~$W_0$ in~$X_0$
    and the hyperplane section~$Z = T \cap H$.
    For simplicity, $Z \cup L$ will be denoted by~$W_1$.
    Then multiplication by $h_{0}$ defines the exact sequence
\begin{eqnarray}
    0 \rightarrow \II_{W_1}(4) \stackrel{\cdot h_{0}}{\longrightarrow}
\II_{W_0 \cup W_1}(5)
    \rightarrow \II_{(W_0 \cup W_1) \cap P, P}(5) \rightarrow 0,
\label{eq:exact_4.0}
\end{eqnarray}
    where~$P = H \cap H_0$.
    The vector space~$H^0 \II_{(W_0 \cup W_1) \cap P, P}(5)$ is
    spanned by the conics in~$P$ containing
    the $4$ points obtained as~$T \cap P - T \cap L$, so $\h^0 \II_{(W_0
    \cup W_1) \cap P, P}(5)=2$.

    On the other hand~$Z \cdot L = 6$,
    so
\[
    \begin{aligned}
    p_a (W_1) = p_a(Z) + p_a(L) + Z \cdot L -1  = 10 - 0 + 6 - 1 = 15.
    \end{aligned}
\]
    By the Riemann-Roch formula,
\[
    \begin{aligned}
    \chi (\OO_{W_1}(4)) = 4 \deg (W_1) + 1 - p_a(W_1) = 44 + 1 - 15 =30,
    \end{aligned}
\]
    and hence~$\chi (\II_{W_1}(4)) =5$,
    which implies that~$\h^0 \II_{W_1} (4) \geq 5$.
    But~$\h^0 \II_{W_1}(4) \leq \h^0 \II_Z (4) \leq 5$,
    so~$\h^0 \II_{W_1}(4) = 5$ and~$\h^1 \II_{W_1}(4)=0$. Thus
    from~(\ref{eq:exact_4.0})
    we have~$\h^0 \II_{W_1 \cup W_2}(5) = 7$.
\end{proof}
\begin{lem}
    The reducible surface~$Y = T \cup X_0$ is scheme-theoretic
    cut out by~$5$ quintic hypersurfaces.
\label{th:global_generated}
\end{lem}
\begin{proof}
    By Lemma \ref{exact}
    the residual sequence
\[
    0 \rightarrow \II_T (4) \stackrel{\cdot h_0}{\longrightarrow}
    \II_Y (5) \rightarrow \II_{D \cup X_0,H_0} (5)
    \rightarrow 0
\]
    remains exact after taking global sections, so the intersection of
    the base
    locus in the middle with $H_{0}$ coincides with the base locus on
    the right. But
    $\II_{D \cup X_0,H_0} (5)$ is globally generated by its own
    two global sections,
    because~$\H^0 \II_{D \cup X_0,H_0} (5) \simeq \H^0 \II_{D,H_0} (2)$,
    so the base locus of the quintics through~$Y$ would be in~$U_0$.
    As we can check in the example given in~\cite{abo}, p.107,
    there is a quintic hypersurface~$V$ through~$Y$ such that
    $V$ does not contain any plane of~$U_0$.
    So this is the general property,
    and hence we may assume the existence of such a quintic
    hypersurface~$V$.
    Note that~$V$ intersects each plane of~$U_0$ in a quartic
    plane curve in~$T$ and the line~$L$ in~$X_0$.
    So~$Y$ is scheme-theoretically cut out by $5$ quintics.
\end{proof}

Now we compute the sectional genus~$\pi$,
the geometric genus~$p_g$ and the irregularity~$q$ of~$X$.
   From the standard liaison sequence
\[
    0 \rightarrow \OO_X(K_X) \rightarrow \OO_{T \cup X_0 \cup X}(5)
    \rightarrow \OO_{T \cup X_0} (5)\rightarrow 0,
\]
we have the exact sequence
\begin{eqnarray}
    0 \rightarrow \II_{T \cup X_0 \cup X}(5) \rightarrow
    \II_{T \cup X_0} (5) \rightarrow \OO_X(K_X) \rightarrow 0.
\label{eq:exact_4.1}
\end{eqnarray}
Twisting by~$\OO(-1)$ and taking cohomology,
we can compute the speciality
\[
     s = \h^1 \OO_X(H_X) = \h^1 \OO_X (K_X - H_X)
       = \h^1 \II_{T \cup X_0}(4) = \h^1 \II_T (3) = 2.
\]
The exact sequence remains exact on global sections,
because~$T \cup X_0 \cup X$ is a complete intersection of two quintics.
Recall, from above,  that
$\h^0 \II_{T \cup X_0}(5) = 5$,
so we have
\[
    p_g = \h^0 \OO_X (K_X) = \h^0 \II_{T \cup X_0} (5)
          - \h^0 \II_{T \cup X_0 \cup X}(5) = 5-2 = 3.
\]
Since~$\chi (\OO_X) = \chi (\OO_X(K_X))
= \chi (\II_{T \cup X_0}(5)) - \chi (\II_{T \cup X_0 \cup X}(5))
= 3$,
we have~$q = 1$.
Now the speciality~$s=\h^1 \OO_X(H_X)=2$ and $\h^2 \OO_X(H_X)=0$, so by
Riemann Roch
$$H_{X}\cdot K_{X}=H_{X}^2+2\chi (\OO_X)-2\chi
(\OO_X(H_{X}))=12+2\cdot 3 -2\cdot (5-2)=12.$$
In particular, by the adjunction formula, $\pi = 13$.
Formula~(\ref{eq:self_intersect}) yields~$K^2 = 0$.

Finally we show that~$X$ is an elliptic surface.
To do this, we show that~$X$ is minimal.
Since~$K^2 = 0$ and~$p_g = 3$,
it is enough to prove that~$K$ has no fixed components.
Let~$C = (X \cup T) \cap X_0$.
By the standard liaison sequence
\[
    0 \rightarrow \OO_{X_0} (K_{X_0}) \rightarrow
    \OO_{X_0}(5) \rightarrow \OO_C(5) \rightarrow 0
\]
and the exact sequence
\[
     0 \rightarrow \OO_{X_0} (-C) \rightarrow
    \OO_{X_0} \rightarrow \OO_C \rightarrow 0,
\]
we can see that~$C \equiv 5H_{X_0} - K_{X_0} \equiv 6H_{X_0}$.
Recall that~$C_1 = T \cap X_0$ is
the residual curve $T \cap H_0 -D$.
So~$C_{1}\equiv 3H_{X_0} - 3L$ on $X_{0}$.  Therefore $C_0 = X \cap
X_{0}\equiv 3H_{X_0} + 3L$.  In particular $C_{0}$  has
degree~$12$, which implies that~$X \cap X_0$ is equal to
the hyperplane section~$H_X = X \cap H_0$.
Applying the analogous argument to $X$ shows that~$(T \cup X_0) \cap
X = 5H_X - K_X$.
Thus $X \cap T \equiv 4H_X - K_X$.
This says that any quartic hypersurface containing~$T$
intersects~$X$ in~$X \cap T$ and~$K_X$.
Note that there is a net of such quartics.
Consider the following commutative diagram
$$
\xymatrix{
& & \II_T(4) \ar[d] \ar[dr] \\
0 \ar[r] & \II_{T\cup X_0 \cup X}(5) \ar[r] & \II_{T\cup X_0}(5) \ar[r]
& \mathcal O_X(K_X) \ar[r] & 0,
}
$$
The global sections of~$\II_T (4)$ are mapped to~$\H^0 \OO_X(K_X)$
injectively, because we have chosen two general quintics
of~$\H^0 \II_{T \cup X_0} (5)$.
Since~$p_g = 3$, we have an isomorphism
\[
    \H^0 \II_T(4) \simeq \H^0 \OO_X (K_X).
\]
This implies that every canonical divisor
of~$X$ appears in the way described above.
Notice that the base locus of the quartics through~$T$
is the union of~$T$ and the three planes~$U_0$.
Therefore any base component of~$|K_X|$ would lie in~$U_0$.
As claimed in the proof of Lemma~\ref{th:global_generated},
any quintic through~$T$ which does not contain
any plane of~$U_0$ intersects each plane of~$U_0$
in a quartic plane curve in~$T$ and the line
in~$X_0$. So if there would be a curve in~$X \cap U_0$,
then this curve would meet~$L$.
On the other hand,
$C_0 \cdot L = (3H_{X_0} + 3L) \cdot L = 0$ on~$X_0$.
This is a contradiction, and hence
$|K_X|$ has no fixed component.
\end{proof}

\begin{rem}
Let~$X$ be an irregular elliptic surface given as above.
We have shown that the canonical system~$|K_X|$ is base point free.
In other words, the canonical bundle
$\omega_X$ is generated by its own three global sections,
because~$p_g~=~3$. This leads to the following:

\noindent
(i) The canonical system~$|K_X|$
defines a morphism~$\Phi :  X \rightarrow \P^2$,
and $\dim \Phi (X) = 1$, because~$K_X^2~=~0$.
Let~$\Phi = s \circ r$ be the Stein factorization of~$\Phi$.
Then, as we saw in the proof of Theorem~\ref{th:monad_const},
the smooth fibers of~$r$ are elliptic curves.

\noindent
(ii) From the generalized Serre correspondence~\cite{okonek}
it follows that there is a rank four vector bundle~$\EE$
with Chern classes~$(c_1,c_2,c_3,c_4) = (5,12,12,0)$
and three sections of~$\EE$ whose dependency locus is~$X$.
However we have not been able to prove that every irregular elliptic surface
we have constructed above arises from a rank three vector bundle
on~$\P^4$ with Chern classes~$(c_1,c_2,c_3) = (5,12,12)$
(compare Remark~\ref{th:rank_three}).
\end{rem}



\begin{rem}
Let~$X$ be a smooth surface in~$\P^4$ with~$d=12$, $\pi=13$,
$p_g=3$ and~$q=1$.  Is it clear that the $X$ is as above?  In this
last remark we sketch an argument why this is true as soon as~$\h^1 
\II_X(2) = 1$.
We conjecture, but cannot prove that this is always the case.
If~$\h^1 \II_X(2) = 1$, then~$\h^2 \II_X(2) = 1$ by the Riemann-Roch theorem.
Using Beilinson's theorem~\cite{beilinson} and  a remark in~\cite{decker},
we can deduce that~$\h^2 \II_X(3)=0$, and thus~$\h^1 \II_X(3)=4$.
This implies that there exists the unique hyperplane~$H_0$
such that~$\h^0 \II_{X \cap H_0}(3) = 1$.
Let~$X_0$ be the unique cubic surface in~$H_0$
containing~$X \cap H_0$.
Then the residual surface~$X \cup X_0$ is locally complete intersection,
because~$X \cap H_0$ and~$X \cap X_0$ coincide.
Suppose that~$X \cup X_0$ is scheme-theoretically cut out by quintics.
Then the linked surface~$T$ is a smooth general type surface
of~$d=10$, $\pi = 10$ and~$\chi =4$.
So~$X$ is reobtained as the surface linked~$(5,5)$
to~$T \cup X_0$, and hence~$X$ is an elliptic surface as above.
\end{rem}


\end{document}